\documentclass[10pt]{amsart}
\usepackage{fix-cm}
\usepackage{times}
\usepackage{amsaddr}
\usepackage{verbatim}
\usepackage{url}
\usepackage{amssymb}
\def\mymedskip{\vskip\medskipamount}
\def\mymedbreak{\par \ifdim\lastskip<\medskipamount
  \removelastskip \penalty-100 \mymedskip \fi}
\def\myaftermedspace{\par \ifdim\lastskip<\medskipamount
  \removelastskip \penalty55\mymedskip\fi}
\newcommand{\eop}{{\unskip\nobreak\hfil\penalty50
          \hskip2em\hbox{}\nobreak\hfil$\Box$
          \parfillskip=0pt \finalhyphendemerits=0 \par}}
{\mymedbreak{\noindent\bf Proof of Theorem #1:\enspace}}{\eop\myaftermedspace}
{\mymedbreak\noindent{\bf Remark:}%
\enspace\rm}%
{\myaftermedspace}
\newtheorem{teor}{Theorem}[section]

\newtheorem{fact}[teor]{Fact}
\newtheorem{problem}{Problem}
\newtheorem{exercise}{Exercise}
\newtheorem{examp}[teor]{Example}
\newtheorem{lem}[teor]{Lemma}
\newtheorem{cor}[teor]{Corollary}

\newtheorem{prop}[teor]{Proposition}
\newtheorem{rem}[teor]{Remark}
\newcommand{\beq}{\begin{equation}}
\newcommand{\eeq}{\end{equation}}
\newcommand{\beql}[1]{\begin{equation} \label{#1}}
\newcommand{\eeql}{\end{equation}}
\newcommand{\beqa}{\begin{eqnarray*}}
\newcommand{\eeqa}{\end{eqnarray*}}
\newcommand{\beqal}[1]{\begin{eqnarray} \label{#1}}
\newcommand{\eeqal}{\end{eqnarray}}
\newcommand{\beqan}{\begin{eqnarray}}
\newcommand{\eeqan}{\end{eqnarray}}
\newcommand{\bpf}{\begin{proof}}
\newcommand{\epf}{\end{proof}}
\newcommand{\ben}{\begin{enumerate}}
\newcommand{\een}{\end{enumerate}}
\newcommand{\bit}{\begin{itemize}}
\newcommand{\eit}{\end{itemize}}

\newcommand{\cE}{{\mathcal E}}

\newcommand{\cQ}{{\mathcal Q}}

\newcommand{\cZ}{{\mathcal Z}}

\newcommand{\gD}{\Delta}
\newcommand{\gre}{\epsilon}
\newcommand{\gl}{\lambda}

\newcommand{\bab}{\begin{abstract}}
\newcommand{\eab}{\end{abstract}}
\newcommand{\bke}{\begin{keywords}}
\newcommand{\eke}{\end{keywords}}

\newcommand{\btm}[1]{\begin{teor} \label{#1}}
\newcommand{\etm}{\end{teor}}
\newcommand{\btmn}[2]{\begin{teor}[#1] \label{#2}}
\newcommand{\etmn}{\end{teor}}
\newcommand{\ble}[1]{\begin{lem} \label{#1}}
\newcommand{\ele}{\end{lem}}
\newcommand{\bLe}[1]{\begin{Lemma} \label{#1}}
\newcommand{\eLe}{\end{Lemma}}
\newcommand{\bpn}[1]{\begin{prop} \label{#1}}
\newcommand{\epn}{\end{prop}}
\newcommand{\bex}[1]{\begin{examp} \label{#1}}
\newcommand{\eex}{\end{examp}}
\newcommand{\bde}[1]{\begin{definition} \label{#1}}
\newcommand{\ede}{\end{definition}}
\newcommand{\bco}[1]{\begin{cor} \label{#1}}
\newcommand{\eco}{\end{cor}}
\newcommand{\bcorn}[2]{\begin{cor}[#1] \label{#1}}
\newcommand{\ecorn}{\end{cor}}
\newcommand{\bcon}[1]{\begin{conjecture} \label{#1}}
\newcommand{\econ}{\end{conjecture}}
\newcommand{\bfa}[1]{\begin{fact} \label{#1}}
\newcommand{\efa}{\end{fact}}
\newcommand{\bpr}[1]{\begin{problem} \label{#1}}
\newcommand{\epr}{\end{problem}}
\newcommand{\bexer}[1]{\begin{exercise} \label{#1}}
\newcommand{\eexer}{\end{exercise}}
\newcommand{\bre}[1]{\begin{rem} \label{#1}}
\newcommand{\ere}{\end{rem}}

\newcommand{\bbQ}{\mathbb{Q}}

\newcommand{\bbZ}{\mathbb{Z}}

\newcounter{question_number}
\newenvironment{question}{\addtocounter{question_number}{1}\noindent{\bf Question \arabic{question_number}}}{\myaftermedspace}

\newenvironment{solution}{\noindent {\bf Solution:} \enspace}{\eop\myaftermedspace}

\newenvironment{hint}{\noindent {\bf Hint:} \enspace}{\eop\myaftermedspace}

\newenvironment{multisolution}[1]{\noindent {\bf Solution #1:} \enspace}{\eop\myaftermedspace}

\newcommand{\bqu}{\begin{question}}
\newcommand{\equ}{\end{question}}
\newcommand{\bs}{\begin{solution}}
\newcommand{\es}{\end{solution}}
\newcommand{\bh}{\begin{hint}}
\newcommand{\eh}{\end{hint}}
\newcommand{\bms}[1]{\begin{multisolution}{#1}}
\newcommand{\ems}{\end{multisolution}}

        \DeclareMathOperator{\GL}{GL}

\newcommand{\te}{\tilde{e}}
\newcommand{\gS}{\Sigma}
\newcommand{\im}{{\rm Im}}

\renewcommand{\ker}{{\rm Ker}}

\begin{document}
%\begin{titlepage}
\title[Critical groups of de Bruijn and Kautz graphs and 
circulants]{Critical groups of generalized de Bruijn and Kautz graphs and 
circulant matrices over finite fields}
\date{\today}
\author{%
Swee Hong Chan, Henk D. L. Hollmann, Dmitrii V. Pasechnik} 
\address{School of Physical and Mathematical Sciences,
Nanyang Technological University, 21 Nanyang Link, Singapore  637371 }
\email{\{sweehong, henk.hollmann, dima\}@ntu.edu.sg}%\\\\
\date{}
\begin{abstract}
We determine the critical groups of the generalized de Bruijn graphs
$\mathrm{DB}(n,d)$ and generalized Kautz graphs $\mathrm{Kautz}(n,d)$, thus
extending and completing earlier results for the classical de Bruijn and Kautz
graphs. Moreover, for a prime $p$ 
the critical groups of $\mathrm{DB}(n,p)$ are shown to be in close
correspondence with groups of $n\times n$ circulant matrices over $\mathbb{F}_p$, 
%HDLH
%explaining numerical data in \cite{oeis:A027362}, and prompting %HDLH (to) ?
%the a possibility
%to construct normal bases in  $\mathbb{F}_{p^n}$ from spanning trees in $\mathrm{DB}(n,p)$.
which explains numerical data in \cite{oeis:A027362}, and suggests the possibility 
to construct normal bases in  $\mathbb{F}_{p^n}$ from spanning trees in $\mathrm{DB}(n,p)$.
\end{abstract}
\maketitle

\vspace{-1cm}
\section{Introduction}
The {\em critical group\/} of a directed graph~$G$ is an abelian group obtained
from the Laplacian matrix $\gD$ of~$G$; it determines and is determined by the Smith
Normal Form (SNF) of~$\gD$. (For precise definitions of these and other terms,
we refer to the next section.) The {\em sandpile group\/} $S(G,v)$ of
$G$  at a vertex $v$ is an abelian group obtained from the 
reduced Laplacian $\gD_v$ of~$G$; its order is equal to the {\em complexity\/}
$\kappa(G)$ of~$G$, the number of directed trees rooted at~$v$, a fact
that is related to the Matrix Tree Theorem, see for example \cite{LeLI-sand}
and its references.  
If $G$ is Eulerian, then $S(G,v)$ does not depend on
$v$, and is then simply written as~$S(G)$; in that case, it is equal to
the critical group of~$G$. The critical group has been studied in other
contexts under several other names, such as  group of components,  Picard or
Jacobian group, and Smith group. For more details and background, see, e.g.,
\cite{HLMPPW-chip}.

Critical groups have been determined for a large number of graph families. For
some examples,  see the references in \cite{AlVa-cone}.  Here, we determine the
critical group of the generalized de Bruijn graphs $\mathrm{DB}(n,d)$ and
generalized Kautz graphs $\mathrm{Kautz}(n,d)$, thus extending and completing
the results from \cite{LeLI-sand} for the binary de Bruijn graphs
$\mathrm{DB}(2^\ell,2)$ and 
Kautz graphs (with $p$ prime) 
$\mathrm{Kautz}((p-1)p^{\ell-1},p)$, and \cite{BiKi-knuth} for the
classical de Bruijn graphs $\mathrm{DB}(d^\ell,d)$ and Kautz graphs
$\mathrm{Kautz}((d-1)d^{\ell-1},d)$. Unlike the classical case, the generalized
versions are not necessarily  iterated line graphs, so to obtain their critical
groups, different techniques have to be applied. 

Our original motivation for studying these groups stems from their relations
to some algebraic objects, such as the groups $C(n,p)$ of invertible
$n\times n$-circulant matrices over
$\mathbb{F}_{p}$ (mysterious numerical coincidences were noted in the OIES
entry A027362 \cite{oeis:A027362} by the third
author, computed with the help of \cite{sage,sagesandpiles}), 
and {\em normal bases} (cf. e.g. \cite{MR1429394}) 
of the finite fields $\mathbb{F}_{p^n}$. The latter were noted
to be closely related to circulant matrices and to {\em necklaces\/} by
Reutenauer \cite[Sect.~7.6.2]{Reu-lie}, see also \cite{DuPa-neck}, and the related
numeric data collected in~\cite{Arn-comp}. 
In particular, we show that 
$C(n,p)/(\mathbb{Z}_{p-1}\times\mathbb{Z}_n)$ is isomorphic to 
the critical group of $\mathrm{DB}(n,p)$. Although we were not able to construct 
an explicit bijection between the former and the latter, we could speculate
that potentially one might be able to design a new deterministic way to construct
normal bases of  $\mathbb{F}_{p^n}$.

\section{Preliminaries}\label{sect:prelims}
Let $M$ be an $m\times n$ integer matrix of rank $r$.
% and nonzero invariant factors $(d_1, \ldots, d_r)$ \cite{Ne-im}.
For a ring $F$, we write $R_F(M)=M^\top F^n$, the $F$-module generated by the rows of $M$. 
% of $\bbQ$, 
%$R_\bbZ(M)$ and $R_\bbQ(M)$ the $\bbZ$-row space and the $\bbQ$-row space of~$M$, respectively, that is, the $\bbZ$-module and the $\bbQ$-module consisting of all linear combinations of the rows of~$M$ with coefficients from~$\bbZ$ and $\bbQ$, respectively. 
The {\em Smith group\/} \cite{Ru-phd}
%, \cite{Ru-evSNF} 
of $M$ is defined as
%HDLH \beql{LESG}
%$
$\Gamma(M)=\bbZ^n/R_\bbZ(M)$.
%$.
%\eeql
%where $R_\bbF(M)$ denotes the $F$-module generated by the rows of $M$. 
The submodule
%\[
$\overline{\Gamma}(M)=\bbZ^n/R_\bbQ(M)\cap \bbZ^n$
%\]
of~$\Gamma(M)$ 
is a finite abelian group called the {\em finite part\/} of $\Gamma(M)$. 
%The justification of this name is that
Indeed, if $M$ has rank $r$, then
%\[
$\Gamma(M) = \bbZ^{n-r} \oplus \overline{\Gamma}(M)$
%\]
with 
%\[
%\[\overline{\Gamma}(M) = \bigoplus_{i=1}^r \bbZ_{d_i},\]
$\overline{\Gamma}(M) = \oplus_{i=1}^r \bbZ_{d_i}$,
where $d_1, \ldots, d_r$ are the nonzero invariant factors of~$M$, so that $d_i|d_{i+1}$ for $i=1, \ldots, r-1$. For invariant factors and the  
%HDLH SNF, 
Smith Normal Form, we refer to~\cite{Ne-im}. See \cite{Ru-phd}
%, \cite{Ru-evSNF} 
for further details and proofs.

Let  $G=(V,E)$ be a directed graph on~$n=|V|$ vertices. The indegree $d^-(v)$ and outdegree 
%HDLH $d^+$ 
$d^+(v)$ is the number of edges ending or starting in $v\in V$, respectively.
%, and let $n=|V||$ denote the number of vertices of~$G$. 
%=\{v_1, \ldots, v_n\}$. 
The {\em adjacency matrix\/} of~$G$ is the $n\times n$ matrix 
%HDLH $A=(A_{v,w}$, 
$A=(A_{v,w})$, with rows and columns indexed by $V$, 
%HDLH Added space after "where"
where $A_{v,w}$ is the number of edges from $v$ to~$w$.
%Let $d^-_v$ denote  the number of edges entering~$v$. 
The {\em Laplacian\/} of~$G$ is the matrix $\gD=D-A$, where $D$ is diagonal with $D_{v,v}=d^-_v$.
%, the {\em indegree\/} of~$v$. the number of edges ending in~$v$. 
The {\em critical group\/} $K(G)$ of~$G$ is the finite part of the Smith group of the Laplacian~$\gD$ of~$G$. 
The {\em sandpile group\/} $S(G,v)$ of~$G$ at a $v\in V$ is the finite part of the Smith group of the $(n-1)\times (n-1)$ {\em reduced Laplacian\/} $\gD_v$, obtained from $\gD$ by deleting the row and the column of $\gD$ indexed by~$v$. Note that by the Matrix Tree Theorem for directed graphs, the order of $S(G,v)$ equals the number of directed spanning trees rooted at~$v$. %HDLH added a "."
%replacing the $v$-th row of $\gD$ by the $v$-th unit vector.
$G$ is called  {\em Eulerian\/} if $d^+(v)=d^-(v)$ for every $v\in V$. 
In that case, $S(G,v)$ does not depend on $v$ and is equal to the critical group $S(G)$ of~$G$. For more details on sandpile groups and critical groups of directed graphs, we refer for example to  \cite{HLMPPW-chip} or \cite{wag-crit}.
% is equal to each of the {\em sandpile groups\/} $S(G,v)$ of~$G$ relative to a vertex~$v$. 
%There are many other ways to define 

\subsection{Generalized de Bruijn and Kautz graphs}
Generalized de Bruijn graphs 
and generalized Kautz graphs  \cite{DCH-dBKgen} are 
known to have a relatively small diameter and attractive connectivity properties, and have been studied intensively due to their applications in interconnection networks. The  generalized Kautz graphs were first investigated  in \cite{ImIt-diam-ieee}, 
%\cite{ImIt-diam-net}, 
and are also known as {\em Imase-Itoh digraphs\/}. Both classes of graphs are Eulerian.

We will determine the critical group, or, equivalently, the sandpile group, of a generalized de Bruijn or  Kautz graph on $n$ vertices by embedding this group as a subgroup of index~$n$ in a group that we will refer to as the {\em sand dune\/} group of the corresponding digraph. 
%This embedding method can in fact be applied to the much wider class of {\em consecutive-$d$ digraphs\/} \cite{DHP-cons}.
Let us now turn to the details.

The {\em generalized de Bruijn graph\/} $\mathrm{DB}(n,d)$ has vertex set $\bbZ_n$, the set of integers modulo $n$, and (directed) edges $v\rightarrow dv+i$ for $i=0, \ldots, d-1$ and all $v\in  \bbZ_n$. The {\em generalized Kautz graph\/} 
%(sometimes also also called {\em Imase-Itoh digraph\/}) 
$\mathrm{Kautz}(n,d)$ has vertex set $\bbZ_n$ and directed edges $v\rightarrow -d(v+1) +i$ for $i=0, \ldots, d-1$ and all $v\in\bbZ_n$. 
%These graphs are important special cases of the so-called {\em consecutive-$d$ digraphs\/}. Here a  consecutive-$d$ digraph $G(d, n, q, r)$, defined for $q\in \bbZ_n\setminus \{0\}$ and $r\in\bbZ_n$, has vertex set $\bbZ_n$ and directed edges $v\rightarrow qv+r +i$ for $i=0, \ldots, d-1$ and all $v\in \bbZ_n$. Note that the generalized de Bruijn and Kautz digraphs are the cases $q=d, r=0$ and $q=-d, r=-d$, respectively. 
Note that both $\mathrm{DB}(n,d)$ and $\mathrm{Kautz}(n,d)$ are Eulerian.
%, and the Imase-Itoh digraph is the case $q=r=-d$. 
%
In what follows, we will 
%HDLH concentrate (and added ";" instead of "."
focus on the generalized de Bruijn  graph;
%HDLH , but we will briefly indicate how to adapt the methods for the  case of the generalized Kautz graph. 
%HDLH Our starting point  will be the Smith group as defined in (\ref{LESG}). 
%HDLH
the generalized Kautz graph can be handled in a similar way, essentially by replacing $d$ by $-d$ in certain places.

Let 
%$\cZ_n =\bbZ[x]\bmod x^n-1$ and $\cQ_n=\bbQ[x]\bmod x^n-1$, and set 
$\cZ_n =\{a(x)\in \bbZ[x]\bmod x^n-1 \mid a(1)=0\}$.
% and $\cQ_n=\{a(x)\in \bbQ[x]\bmod x^n-1 \mid a(1)=1\}$. 
With each vertex $v\in \bbZ_n$, we associate the  polynomial $f_v(x) =  dx^v -x^{dv}\sum_{i=0}^{d-1} x^{i}\in\cZ_n$.
Since $f_v(x)$ is the associated polynomial of the $v$th row of the Laplacian $\gD^{(n,d)}$ of the generalized de Bruijn graph $\mathrm{DB}(n,d)$, 
%HDLH we see from (\ref{LESG})  that
the Smith group $\Gamma(\Delta^{(n,d)})$ of the Laplacian of $\mathrm{DB}(n,d)$ is the quotient 
%\beql{LESGDB}
%\Gamma(\Delta^{(n,d)}) = \cZ_n/\langle f_v(x) \mid v\in \bbZ_n\rangle_{\bbZ_n} 
%\Gamma(\Delta^{(n,d)}) = (\bbZ[x] \bmod x^n-1)/\langle f_v(x) \mid v\in \bbZ_n\rangle_{\bbZ_n} 
%\eeql
%of $\cZ$ 
of $\bbZ[x] \bmod x^n-1$
by the $\bbZ_n$-span $\langle f_v(x) \mid v\in \bbZ_n\rangle_{\bbZ_n}$ of the polynomials $f_v(x)$. Now note that 
%\[\cZ_n =\bbZ \oplus \cZ_n',\]
%\[
$\bbZ[x] \bmod x^n-1 \cong \bbZ \oplus \cZ_n$,
%\]
%where $\cZ_n'=\{a(x)\in \cZ_n \mid a(1)=1\}$. 
so since 
%$f_v(1)=0$ and 
$\sum_{v\in \bbZ_n} f_v(x)=0$, we have that
%\[  \cZ_n/\langle f_v(x) \mid v\in \bbZ_n\rangle_{\bbZ_n} \cong \bbZ \oplus \cZ_n'/ \langle f_v(x) \mid v\in \bbZ_n'\rangle_{\bbZ_n}\]
\begin{align}
 \Gamma(\Delta^{(n,d)})
&= (\bbZ[x] \bmod x^n-1)/\langle f_v(x) \mid v\in \bbZ_n\rangle_{\bbZ_n}  \\& \cong \bbZ \oplus \cZ_n/ \langle f_v(x) \mid v\in \bbZ_n'\rangle_{\bbZ_n} \notag
\end{align}
where $\bbZ_n'=\bbZ_n\setminus \{0\}$. It is easily checked that the polynomials $f_v(x)$ with $v\in \bbZ_n'$ are independent over $\bbQ$, hence they constitute a basis for $\cZ_n$ over $\bbQ$. As a consequence, each element in the quotient group
\beql{LESPDB}S(n,d)=  S_{\mathrm{DB}}(n,d) = \cZ_n/ \langle f_v(x) \mid v\in \bbZ_n'\rangle_{\bbZ_n}\eeql
%HDLH1 (added "equivalently")
has finite order, and so $S(n,d)$ is the critical group, or, equivalently,  the sandpile group, 
of the generalized de Bruijn graph  $\mathrm{DB}(n,d)$.
We define the {\em sand dune group\/} $\gS(n,d)=\gS_{\mathrm{DB}}(n,d)$ of $\mathrm{DB}(n,d)$ 
%is defined
as
% in terms of the polynomials 
%\[
%$ g_v(x) = (x-1)f_v(x) = dx^v (x-1) -x^{dv}(x^d-1)$
%\]
%\[ 
$\gS(n,d) =  
%\gS_{\mathrm{DB}}(n,d) =  
\cZ_n/ \langle g_v(x) \mid v\in \bbZ_n'\rangle_{\bbZ_n}$,
%\]
where $ g_v(x) = (x-1)f_v(x) = dx^v (x-1) -x^{dv}(x^d-1)$.
%HDLH1 Below, I have replaced all \Env{n}{v} (e^n) by just e and all Grenv{n}{v} (epsilon^(n)) by just epsilon...
%begin %HDLH1
%
Now let
%$\Env{n}{v}=x^v-1$; 
$e_v=x^v-1$; 
we have that 
% $\Env{n}{0}=0$, 
$e_0=0$, and $\cZ_n=\langle 
%\Env{n}{v} 
e_v \mid v \in \bbZ_n'\rangle_{\bbZ}$, the $\bbZ$-span of the
%$\Env{n}{v}$. 
polynomials~$e_v$. 
%HDLH1 For later use, we also define (replace this:
Furthermore, let 
%\[ 
%$\Grenv{n}{v}=d\Env{n}{v}-\Env{n}{dv}$.
$\gre_v=de_v-e_{dv}$.
%\]
%In what follows, we will simply write $e_v$ and $\gre_v$ instead of $\Env{n}{v}$ and $\Grenv{n}{v}$ if the value of~$n$ %is evident from the context.
%
The span in $\cQ_n=\{a(x)\in \bbQ[x]\bmod x^n-1 \mid a(1)=0\}$ of the polynomials $g_v(x)$  with $v\in \bbZ_n$ is the set of polynomials of the form $dc(x)-c(x^d)$ with $c(1)=0$; since 
%\[
$\gre_v = g_0(x) + \cdots g_{v-1}(x)$
% \]
for all $v\in \bbZ_n$, we conclude that
\beql{LEsd-def}\gS(n,d)=\cZ_n/\cE_{n,d},\eeql
%[note the change in notation, no Z' but Z]
where 
%HDLH1 added for more clarity
 $\cZ_n=\langle e_v \mid v \in \bbZ_n'\rangle_{\bbZ}$ and
%\[  
$\cE_{n,d}=\langle 
%d\Env{n}{v}-\Env{n}{dv}
\gre_v\mid v\in \bbZ_n'\rangle_\bbZ$
%\]
is the $\bbZ$-submodule of $\cZ_n$ generated by the polynomials 
%\[ %\Endvpol{n}{d}{v}
%$\Grenv{n}{v}$.
%HDLH1 added for clarity
$\gre_v=de_v-e_{dv}$.
%=\Endv{n}{d}{v}.\]
%
% %End %HDLH1
%
The next result is crucial: it identifies the elements of 
the sand dune group
$\gS(n,d)$ that are actually contained 
in the sandpile group $S(n,d)$.  (Due to lack of space, we omit  the not too difficult 
 proofs %, which are not difficult, 
 in the remainder of this section.)
\btm{Lsub}
If $a\in \gS(n,d)$ with $a=\sum_v a_v e_v$, then $a\in S(n,d)$ if and only if $\sum_v va_v\equiv0 \bmod n$.
\etm
%remove the proof of theorem 2.1
\iffalse
\bpf
%Let $a\in S(n,d)$. 
The 
%``multiplication by $(x-1)$'' 
map $\phi: a(x) \rightarrow (x-1)a(x)$ on $\cZ_n=\langle e_{v} \mid v \in \bbZ_n'\rangle_{\bbZ}$ 
%that sends  $a(x)$ to  $(x-1)a(x)$ 
is a homomorphism on $(\cZ_n,+)$ mapping $ \langle f_v(x) \mid v\in \bbZ_n\rangle_{\bbZ_n}$ one-to-one onto~$\cE_{n,d}$, hence it embeds $S(n,d)$ as a subgroup $\phi(\cZ_n)/\cE_{n,d}$ in $\gS(n,d)$. To determine that subgroup, we need to determine the image of $\phi$ on $\cZ_n$.

Let $a(x)\in \cZ_n$. Then $a(1)=0$, so $a(x)$ is divisible by $x-1$. Write $a'(x)=a(x)/(x-1)$. Note that $a'(x)$ has integer coefficients, and that  $a(x)=b(x)(x-1)$ in $\bbQ[x]\bmod x^n-1$ iff $b(x)=a'(x)-\gl T(x)$ with $\gl\in \bbQ$, where $T(x)=(x^n-1)/(x-1)$.   Now $T(1)=n$, hence $b(1)=0$ iff $\gl= a'(1)/n$, and so $b(x)\in\cZ_n$ iff $\gl$ is integer, that is, iff $a'(1)\equiv 0\bmod n$. 
To finish the proof, note that $T_v(x)=e_v(x)/(x-1)=1+x+\cdots +x^{v-1}$, hence $T_v(1)=v$; now
$a'(x)=\sum_v a_vT_v(x)$ so $a'(1)=\sum_vva_v$.
\epf
\fi
\bco{LCord}
We have  $\gS(n,d)/S(n,d)=\bbZ_n$ and so
%\[
$|\gS(n,d)|=n|S(n,d)|$.
%\]
\eco
%CSH Removing everything from Corollary 2.2 onwards
%HS\iffalse

\iffalse %HS
\bpf

The map $\phi:\sum_v a_v e_v \rightarrow \sum_vva_v \bmod n$ has the property that $\phi(\gre_v)=\phi(de_v-e_{dv})=0$, hence is well-defined as a map on~$\gS(n,d)$. The map $\phi$  is obviously a homomorphism, and is  surjective  since $\phi(e_v)=v$ for all $v\in \bbZ_n$. 
As a consequence, $n=|\bbZ_n|=|\im(\phi)|=|\gS(n,d)|/|\ker(\phi)|=|\gS(n,d)|/|S(n,d)|$.
\epf
\fi %HS
%HDLH1 extra blabla
The above descriptions of the sandpile group $S(n,d)$ and sand dune group $\gS(n,d)$, and the embedding of $S(n,d)$ as a subgroup of $\gS(n,d)$ are very suitable for the determination of these groups. In the process, repeatedly information is required about the order of various group elements. The following two results provide that information.
%
%HDLH1 removed now
%The following result is useful to determine the {\em order\/} of elements in the embedding group $\gS(n,d)$. 
\ble{Lord} Let $a=\sum_v a_v\gre_v\in \gS(n,d)$. Then the order of $a$ in $\gS(n,d)$ is the smallest positive integer~$m$ for which $ma_v\in \bbZ$ for each~$v$.
\ele
\iffalse %HS
\bpf
Since the $\gre_v$ are independent in $\bbQ[x]\bmod x^n-1$, each polynomial $f(x)$ in $\bbQ[x]\bmod x^n-1$ with $f(1)=0$ has a unique expression $f=\sum_vf_v\gre_v$ as linear combination of the $\gre_v$. Such an expression is 0 modulo $\cE(n,d)$ iff all coefficients $f_v$ are integers. Now the claim is obvious.
\epf
\fi %HS
%By expressing $e_v$ in terms of the $\gre_v$, 
We say that $v\in \bbZ_n$ has {\em $d$-type\/} $(f,e)$ in~$\bbZ_n$  if $v,dv,\ldots,d^{e+f-1}v$ are all distinct, with $d^{e+f}v=d^fv$.
Now, by expressing $e_v$ in terms of the $\gre_v$, we can determine the order of $e_v$. The result is as follows.
%As we shall see, the order of~$e_v$ is completely determined by the type of~$v$.
\ble{Lexp} Supposing $v$ has $d$-type $(f,e)$,
% in $\bbZ_n$ 
then 
%\[e_v=d^{-1}\gre_v+\cdots+d^{-f}\gre_{d^{f-1}v} + \frac{1}{d^{f+e-1}(d^e-1)}\gre_{d^{f+e}v}+ \cdots + \frac{1}{d^{f}(d^e-1)}\gre_{d^{f+e}v}. \]
%\[
$e_v =\sum_{i=0}^{f-1} d^{-i-1}\gre_{d^iv}+\sum_{j=0}^{e-1}d^{j-f}(d^e-1)^{-1}\gre_{d^{f+j}v}$
in~$\cZ_n$, and hence $e_v$ has order $d^f(d^e-1)$ in $\gS(n,d)$.
%\]
\ele 

%HS\fi

\subsection{Invertible circulant matrices}\label{Icm}
Let $Q_n$ be the $n \times n$ permutation matrix over a field $F$  
corresponding to the cyclic permutation $(1,2,\dots,n)$%
%HS
\iffalse 
, that is $Q_n = \left(\begin{smallmatrix}
0 & 1 &  \ldots& 0 & 0 \\ 
 0 & 0 & 1 &  & 0 \\ 
 \vdots & 0 & 0  &\ddots   &\vdots \\ 
 0&  & \ddots  & \ddots & 1 \\ 
 1& 0 &\ldots  & 0  &0 \end{smallmatrix}\right)$ 
\fi %HS
 %CSH1 remove the matrix to save space 
.
An $n \times n$ \emph{circulant matrix} over $F$ is a matrix that can be
written as $a_1Q_n+a_2 Q_n^2+ \ldots + a_n Q_n^n$ with $a_i \in F$ for $1\leq i\leq n$.
All the invertible circulant matrices form a commutative group (w.r.t. matrix
multiplication), namely, the
centralizer of $Q_n$ in $\GL_n(F)$. In the case $F=\mathbb{F}_p$ we
consider here we denote this commutative group by $C(n,p)$. Note that $C(n,p)$
contains a subgroup isomorphic to $\mathbb{Z}_{p-1}\oplus
\mathbb{Z}_n$, namely the direct product of the group of scalar matrices 
$F_p^* I:=\{\lambda I\mid \lambda\in\mathbb{F}_p^*\}$ and
the cyclic subgroup generated by $Q_n$. Each circulant matrix has all-ones vector
$\mathbf{1}:=(1,\dots,1)^\top$ as an eigenvector. Thus $C'(n,p):=\{g\in C(n,p)\mid 
g\mathbf{1}=\mathbf{1}\}$ is a subgroup of $C(n,p)$, and we have the following formula.
\begin{equation}\label{Cnpdec}
C(n,p)=C'(n,p)\times F_p^* I.
\end{equation}

\section{Main results}
Let $n,  d>0$ be fixed integers. The description of the sandpile group
$S(n,d)$ and the sand-dune group $\gS(n,d)$ of the generalized the Bruin graph
$\mathrm{DB}(n,d)$ involves a sequence of numbers defined as follows.  Put
$n_0=n$, and for $i=1, 2,\ldots$, define $g_i= \gcd(n_i,d)$ and $n_{i+1}=n_i/g_i$.
We have $n_0>\cdots>n_k=n_{k+1}$, where $k$ is the smallest integer for
which $g_k=1$.  
We will refer to the sequence $n_0>\cdots >n_k=n_{k+1}$  as the {\em $d$-sequence\/} of $n$.  
In what follows, we will write $m=n_k$ and
$g=g_0\cdots g_{k-1}$. Note that $n=g m$ with $\gcd(m,d)=1$.

Since $\gcd(m,d)=1$, the map $x\rightarrow dx$ partitions $\bbZ_m$ into orbits  of the form $O(v)=(v,dv,\ldots, d^{o(v)-1}v)$. 
% 
%CSH remove Here, $O(v)$ is sometimes referred to as the {\em $d$-cyclotomic coset of~$v$ modulo $m$\/}.  
We will refer to $o(v)=|O(v)|$ as the 
%{\em order of the $d$-orbit of $v$ on $\bbZ_m$\/}, or briefly, the 
{\em order\/} of $v$. 
%The $d$-orbits are sometimes referred to as {\em $d$-cyclotomic cosets modulo $m$\/}.

For every prime $p|m$, we define $\pi_p(m)$ to be the largest power of $p$ dividing $m$. Let $V$ be a complete set of representatives of the orbits $O(v)$ different from $\{0\}$, where we ensure that for every divisor $p$ of $m$, all integers 
%$p^i m/\pi_p(m)<m$ 
of the form $m/p^j$ 
are contained in~$V$.
% for all $i\geq0$ for which $1\leq p^i<\pi_p(m)$. 
%CSH remove to save space(This is possible since no two of these numbers are in the same coset.) 

\btm{LTmain} With the above definitions and notation, we have that
\beql{LEdb-sd} \gS(n,d)= \biggl[ \bigoplus_{i=0}^{k-1} \bbZ_{d^{i+1}}^{n_i-2n_{i+1}+n_{i+2}}\biggl] \oplus \biggl[ \bigoplus_{v\in V}\bbZ_{d^{o(v)}-1}  \biggr],
\eeql 
 and 
\beql{LEdb-sp}S(n,d)=
\biggl[ \bigoplus_{i=0}^{k-1} \bbZ_{d^{i+1}/g_i}\oplus \bbZ_{d^{i+1}}^{n_i-2n_{i+1}+n_{i+2}-1}\biggl] \oplus 
%\left[\bigoplus_{p|n} \bbZ_{(d^{o(n/\pi_p(n))}-1)/\pi_p(n)}\right]  
%\bigoplus 
\left[  \bigoplus_{v\in V} \bbZ_{(d^{o(v)}-1)/c(v)}\right],
\eeql
where $c(v)=1$ except in the following cases. For any $p|m$,
\[
c(m/\pi_p(m))=
\left\{
\begin{array}{ll}
\pi_p(m), & \mbox{if $p\neq 2$ or $d \equiv 1 \bmod 4$ or $4\!\!\not
|\,m$};\\
\pi_2(m)/2, & \mbox{if $p=2$ and $d\equiv 3 \bmod 4$ and $4|m$},
\end{array}
\right.
\]
and if $4|m$ and $d\equiv 3\bmod 4$, then $c(m/2)=2$.
\etm
%HDLH End replacings
%CSH remove the explanation to save space
\begin{comment}
We remark that since $n=g_0\cdots g_{k-1} m$ with $m=\prod_{p|m} \pi_p(m)$, the above result implies that $ \gS(n,d)/S(n,d)=\bbZ_n$, which can also be seen directly from the results in 
Section~\ref{sect:prelims}.
\end{comment}

For the generalized Kautz graph, a similar result holds. For $v\in \bbZ_m$, we let $O'(v)$ denote the orbit of $v$ under the map $x\rightarrow -dv$, and we define $o'(v)=|O'(v)|$. Now take $V'$ to be a complete set of representatives of the orbits on $\bbZ_m'$. Finally, define $c'(v)$ similar to $c(v)$, except that now $d$ is replaced by $-d$ (so the special case now involves $d\equiv 1\bmod 4$). Then we have the following.
\btm{LTmain-kautz}
The sandpile group $S_{\rm Kautz}(n,d)$ of the generalized Kautz graph
$\mathrm{Kautz}(n,d)$ is obtained from $S(n,d)$ by  replacing  $V$ by $V'$,
$o(v)$ by $o'(v)$, and $c(v)$ by $c'(v)$ in~(\ref{LEdb-sp}).
\etm

%\begin{comment}
The above results can be proved in a number of steps. In what follows, we
outline the method for the generalized de Bruijn graphs; for the generalized
Kautz graphs, a similar approach can be used. Furthermore, we note that many of
the  steps below repeatedly use Theorem~\ref{Lsub} and Lemma~\ref{Lexp}. First,
we investigate  the ``multiplication-by-$d$'' map $d: x\rightarrow dx$ on the
sandpile and sand-dune group. Let $\gS_0(n,d)$ and $S_0(n,d)$ denote the kernel
of the map $d^k$ on $\gS(n,d)$ and $S(n,d)$, respectively.  It is not difficult
to see that $\gS(n,d)\cong \gS_0(n,d)\oplus \gS(m,d)$ and $S(n,d)\cong S_0(n,d)
\oplus S(m,d)$. Then, we use the map $d$ to determine $\gS_0(n,d)$ and
$S_0(n,d)$. It is easy to see that for {\em any\/} $n$, we have $d\gS(n,d)\cong
\gS(n/(n,d), d)$ and $dS(n,d)\cong S(n/(n,d), d)$. With much more effort, it
can be show that the kernel of the map $d$ on $\gS(n,d)$ and $S(n,d)$ is
isomorphic to $\bbZ_d^{n-n/(n,d)}$ and $\bbZ_{d/(n,d)}\oplus
%HDLH Added "-1" in the exponent
\bbZ_d^{n-1-n/(n,d)}$, respectively. Then we use induction over the length $k+1$
of the $d$-sequence of~$n$ to show that $\gS_0(n,d)$ and $S_0(n,d)$ have the
form of the left part of the right hand side in (\ref{LEdb-sd}) and (\ref{LEdb-sp}),
respectively. This part of the proof, although much more complicated, resembles
the method used by \cite{LeLI-sand} and \cite{BiKi-knuth}. 

Now it remains to handle the parts $\gS(m,d)$ and $S(m,d)$ with $\gcd(m,d)=1$.
For the ``helper'' group $\gS(m,d)$ that embeds $S(m,d)$, this is trivial: it
is easily seen that $\gS(m,d)= \oplus_{v\in V} \langle e_v\rangle$, and the
order of $e_v$ is equal to the size $o(v)$ of its orbit $O(v)$ under the map
$d$, so  (\ref{LEdb-sd}) follows immediately. The $e_v$ are not contained
in~$S(m,d)$, but we can try to modify them slightly to obtain a similar
decomposition for $S(m,d)$. The idea is to replace $e_v$ by a modified version
$\te_v=e_v-\sum_{p|m} \gl_p(v) e_{\pi_p(v)m/\pi_p(m)}$, where the numbers $\gl_p(v)$
are chosen such that $\te_v\in S(m,d)$, or by a suitable multiple of $e_v$, in
some exceptional cases (these are cases where $c(v)>1$). It turns out that
this is indeed possible, and in this way the proof of Theorem~\ref{LTmain} can
be completed.

Finally, with the notation from Subsect.~\ref{Icm}, we have the following isomorphisms,
connecting critical groups and circulant matrices.
\btm{LTmain-circ} Let $d$ be a prime. Then %$S(n,d)\oplus\mathbb{Z}_{d-1}\cong C(n,d)/\langle Q_n\rangle $ and
\begin{align} 
% CSH remove this line 
%S(n,d)\oplus\mathbb{Z}_{d-1}\cong C(n,d)/\langle Q_n\rangle , \quad 
 S(n,d)\cong
C'(n,d)/\langle Q_n\rangle  ,  \qquad\text{and} 
 %HDLH added this one:
 \qquad  \gS(n,d)\cong C'(n,d). \notag 
\end{align}
\etm

The proof of Theorem~\ref{LTmain-circ} is by reducing to the case $\gcd(n,p)=1$ by 
an explicit construction, and then by diagonalizing 
$C(n,p)$ over an appropriate extension of $\mathbb{F}_p$. Essentially, as soon as
$\gcd(n,p)=1$, one can read off a decomposition of $C(n,p)$ into cyclic factors from the 
irreducible factors of the polynomial $x^n-1$ over $\mathbb{F}_p$.

%\end{comment}
\bibliographystyle{abbrv} 
\bibliography{sandpile-euro}
\end{document}